\documentclass[11pt]{article}

\usepackage{amsthm, amsmath, amssymb, amsfonts, url}
\usepackage[margin = 1in]{geometry}
\usepackage{tikz}

\usepackage{setspace}
\usepackage{hyperref}

\hypersetup{
    pdfstartview={FitH}    % fits the width of the page to the window
}

\newtheorem*{claim*}{Claim}
\newtheorem*{fact*}{Fact}
\newtheorem{theorem}{Theorem}[section]
\newtheorem*{theorem*}{Theorem}
\newtheorem{proposition}[theorem]{Proposition}
\newtheorem{lemma}[theorem]{Lemma}
\newtheorem*{lemma*}{Lemma}
\newtheorem{corollary}[theorem]{Corollary}

\newtheorem{observation}[theorem]{Observation}

\theoremstyle{definition}

\newtheorem*{question*}{Question}

\theoremstyle{remark}

\def\c{\mathcal}

\def\ZZ{\mathbb{Z}}

\def\Zm#1{\ZZ / #1 \ZZ}

\DeclareMathOperator{\MSTD}{MSTD}
\DeclareMathOperator{\ord}{ord}
\DeclareMathOperator{\trace}{trace}

\newcommand{\abs}[1]{\left\lvert#1\right\rvert}

\newcommand{\sq}{\ \square \ }
\def\({\left(}
\def\){\right)}

\def\x{\times}
\def\vphi{\varphi}

\def\<={\Leftarrow}
\def\=>{\Rightarrow}

\begin{document}

\title{Counting MSTD Sets in Finite Abelian Groups}
\author{Yufei Zhao \\[6pt]
Massachusetts Institute of Technology \\
\url{yufeiz@mit.edu}}

\maketitle

\begin{abstract}
In an abelian group $G$, a more sums than differences (MSTD) set is a subset $A \subset G$ such that $\abs{A+A} > \abs{A-A}$. We provide asymptotics for the number of MSTD sets in finite abelian groups, extending previous results of Nathanson. The proof contains an application of a recently resolved conjecture of Alon and Kahn on the number of independent sets in a regular graph.
\end{abstract}

%\begin{keyword}
%Sumset \sep difference set \sep MSTD \sep finite abelian group \sep independent set
%\MSC[2000] 11P99 \sep 05C69
%\end{keyword}

\tikzstyle{P} = [draw, circle, black, fill, inner sep = 0pt, minimum width = 3pt]
\tikzstyle{every loop} = []

%%%%%%%%%%%%%%%%%%%%%%%%%%%%%%%%%%%%%%%%%%%%%%%%%%%%%%%%%%%%%%%%%%%%%%%%%%%%%%%%%%%%%%%%%%%%%%%%

\section{Introduction} \label{sec:intro}

Let $G$ be an abelian group. A \emph{more sums than differences} set is a finite subset $A \subset G$ such that $\abs{A+A} > \abs{A-A}$, where the sum set $A+A$ and the difference set $A-A$ are subsets of $G$ defined as
\begin{align*}
	A + A &= \{a + b : a, b \in A \}, \\
	A - A &= \{a - b : a, b \in A \}.
\end{align*}
We will use $\MSTD(G)$ to denote the collection of all MSTD sets in $G$, and $\abs{\MSTD(G)}$ the number of MSTD sets in $G$. Since addition is commutative while subtraction is not, two generic elements generate one sum but two differences. So we should generally expect there to be at least as many differences as sums. However, this is not always the case, and the exceptions are the MSTD sets.

There has been extensive of research on MSTD sets of integers. The first example of an MSTD set of integers was found by Conway in the 1960's: $\{0,2,3,4,7,11,12,14\}$. The name MSTD was later given by Nathanson \cite{Na}. For recent papers on MSTD sets, see \cite{Hegarty, HM, MO, MOS, Na2, Na, Zhao:bidirectional, Zhao:limit}. For older papers see \cite{HRY, Marica, PF, Roesler, Ru1, Ru2, Ru3}. We refer the readers to \cite{Na2, Na} for the history of the problem.

Almost all previous research on MSTD sets focused exclusively on the integers, as opposed to other abelian groups. The only paper where MSTD sets in finite abelian groups are considered is by Nathanson \cite{Na}, who showed that families of MSTD sets of integers can be constructed from MSTD sets in finite abelian groups. As an illustration, if $A \subset \Zm n$ satisfies $\abs{A + A} > \abs{A - A}$, then $\{ a \in \ZZ : a \bmod n \in A, 0 \leq a < kn \} $ is an MSTD set of integers for sufficiently large $k$. In general, starting from an MSTD set in a finite abelian group, we can produce an MSTD subset of the lattice $\ZZ^d$ in an analogous manner, and then we can obtain an MSTD set of the integers through the group homomorphism $\phi : \ZZ^d \to \ZZ$ defined by
\[
	\psi(a) = \sum_{i=1}^d a_i m^{i-1}
\]
where $m$ is some sufficiently large integer. This connection to MSTD sets of integers is the initial motivation for studying MSTD sets in finite abelian groups. Nathanson showed that
\[
	\abs{\MSTD(\Zm n \times \Zm 2 )}
	\geq
	\begin{cases}
		2^n \(1 - \frac{2 n}{2^{n/2}} \) & \text{if } n \text{ is even}, \\
	 	2^n \(1 - \frac{\sqrt 2 n}{2^{n/2}} \) & \text{if } n \text{ is odd}.	 
	\end{cases}
\]
Our main result improves and generalizes Nathanson's result. We give asymptotics for $\abs{\MSTD(G)}$ for large finite abelian groups $G$. Recall that the notation $f_n \sim g_n$ means that $f_n / g_n \to 1$ as $n \to \infty$.

\begin{theorem} \label{thm:main}
Let $\{G_n\}$ be a sequence of finite abelian groups with $\abs{G_n} \to \infty$.
\begin{enumerate}
	\item (Even case) If $G_n$ has $k_n > 0$ elements of order 2, and $\limsup_{n \to \infty} \frac{k_n}{\abs{G_n}} < 1 - \frac{1}{2}\log_3 7 = 0.114\dots$, then
		\[
			\abs{\MSTD(G_n)} \sim k_n \cdot 3^{\abs{G_n}/2}.
		\]
	\item (Odd case) If every $\abs{G_n}$ is odd, and the proportion of elements in $G_n$ with order less than $\log_\vphi \abs{G_n}$ approaches $0$ as $n \to \infty$, then
		\[
			\abs{\MSTD(G_n)} \sim \frac 1 2 \abs{G_n} \vphi^{\abs{G_n}} ,
		\]
		where $\vphi = \frac{1 + \sqrt 5}{2}$ is the golden ratio.
\end{enumerate}
\end{theorem}

The hypotheses on $G_n$ in both cases are rather mild. For instance, if the rank of $G_n$ is uniformly bounded, then $k_n$ is uniformly bounded so that the hypotheses for the even case hold, and Proposition~\ref{prop:rank} will show that the hypotheses for the odd case hold as well. As examples, here are two corollaries, the second of which is an improvement of Nathanson's result stated before Theorem~\ref{thm:main}.

\begin{corollary}
$\displaystyle \abs{\MSTD(\Zm n)} \sim 
	\begin{cases}
	 	3^{n/2}, & \text{if } n \text{ is even}, \\
	 	\frac 1 2 n \vphi^n,  & \text{if } n \text{ is odd}.
	\end{cases}$
\end{corollary}

\begin{corollary}
$\displaystyle \abs{\MSTD(\Zm n \times \Zm 2)} \sim 
	\begin{cases}
		3^{n+1} & \text{if } n \text{ is even}, \\
	 	3^{n} & \text{if } n \text{ is odd}.
	\end{cases}$
\end{corollary}

%\begin{corollary} For a fixed positive integer $r$,
%$\displaystyle \abs{\MSTD((\Zm n)^r)} \sim 
%	\begin{cases}
%	 	(2^r - 1) 3^{n} & \text{if } n \text{ is even}, \\
%	 	\frac 1 2 n^r \vphi^{n^r} & \text{if } n \text{ is odd}.
%	\end{cases}$
%\end{corollary}

Let us mention some analogous results for MSTD sets of integers. Martin and O'Bryant \cite{MO} showed that, as $n \to \infty$, there is a positive lower bound to the proportion of subsets of $\{1, 2, \dots, n\}$ which are MSTD sets. Recently the author \cite{Zhao:limit} showed that this proportion in fact approaches a limit greater than $4 \x 10^{-4}$ as $n \to \infty$. Monte Carlo experiments suggest that the limit should be around $4.5 \x 10^{-4}$. Families of MSTD subsets of $\{1, 2, \dots, n\}$ were constructed by Miller, Orosz, and Scheinerman \cite{MOS} and the author \cite{Zhao:bidirectional}.

The rest of the paper contains the proof of Theorem \ref{thm:main}. One feature of the proof is that it contains an application of a graph theory conjecture of Alon \cite{Alon91} and Kahn \cite{Kahn} which was recently resolved by the author \cite{Zhao:indep}; see Theorem \ref{thm:indep}.

The paper is structured as follows. In Section \ref{sec:union-bound} we describe our strategy for bounding the number of MSTD sets. We reduce the problem into counting the number subsets of $G$ that avoid specific sums and/or differences, and they are analyzed using forbiddance graphs, which are discussed in Section \ref{sec:forbiddance}. In Section \ref{sec:bounding} we put all the individual forbiddance graph inequalities together to give a bound for the number of MSTD sets. Section \ref{sec:discussion} contains some further questions.

%%%%%%%%%%%%%%%%%%%%%%%%%%%%%%%%%%%%%%%%%%%%%%%%%%%%%%%%%%%%%%%%%%%%%%%%%%%%%%%%%%%%%%%%%%%%%%%

\section{Union bounds} \label{sec:union-bound}

First, let us give some intuition on why MSTD sets in finite abelian groups behave differently from MSTD sets of integers. In the case of integers, say if we were interested in subsets of $\{1, 2, \dots, n\}$ that are MSTD sets, the ``middle'' sums and differences are almost always present since they can each be represented as a sum or difference in many ways. Thus it is the fringe elements that matter the most. It was recently shown \cite{Zhao:limit} that a ``typical'' MSTD subset of $\{1, 2, \dots, n\}$ consists of a well-controlled fringe and an almost unrestricted middle. On the other hand, in the case of finite abelian groups, there are no longer any fringe elements. Consequently, almost all subsets $A$ of a finite abelian group $G$ satisfy $A + A = A - A = G$, and MSTD sets occupy a diminishing fraction of the subsets of $G$, unlike the integers case.

Now, returning to the case of finite abelian groups, if $A \subset G$ is an MSTD set, then we must have $A - A \neq G$. On the other hand, if $A - A \neq G$ but $A + A = G$, then $A$ is necessarily an MSTD set. Therefore
\begin{equation} \label{eq:containment}
	\{A \subset G : A - A \neq G, A + A = G \} \subset \MSTD (G) \subset \{ A \subset G : A - A \neq G \}.
\end{equation}
We will use union bounds to estimate the sizes of the sets in \eqref{eq:containment}. Let $G'$ denote a subset of $G$ such that for every nonzero $d \in G$, exactly one of $d, -d$ appears in $G'$. Also assume that $0 \notin G'$. For example, if $G = \ZZ/8\ZZ$, then we could use $G' = \{1, 2, 3, 4\} \subset G$ (though $G' = \{1, 3, 4, 6\}$ is an equally valid choice).

Using \eqref{eq:containment}, we obtain the following upper bound
\begin{equation} \label{eq:upper-union}
	\abs{\MSTD (G)} 
	\leq \abs{ \{ A \subset G : A - A \neq G \}}
	\leq \sum_{d \in G'} \abs{\{ A \subset G : d \notin A - A \}}.
\end{equation}

For $d \in G$, let
\[
	D_d = \{ A \subset G : d \notin A - A, A + A = G \}.
\]
Then
\[
	\{A \subset G : A - A \neq G, A + A = G \} = \bigcup_{d \in G'} D_d.
\]
We have
\[
	\abs{\bigcup_{d \in G'} D_d} 
	\geq \sum_{d \in G'} \abs{D_d} - \sum_{\substack{d_1, d_2 \in G' \\ d_1 \neq d_2}} \abs{D_{d_1} \cap D_{d_2}}
\]
Also,
\begin{align*}
	\abs{D_d} &= \{ A \subset G : d \notin A - A, A + A = G \}
		\\  &\geq \abs{ \{ A \subset G : d \notin A - A \}  } 
					- \sum_{s \in G}  \abs{ \{ A \subset G : d \notin A - A, s \notin A + A \}  }, 
\end{align*}
and
\[
	\abs{D_{d_1} \cap D_{d_2}} 
	= \abs{ \{ A \subset G : d_1, d_2 \notin A - A, A + A = G \}  }
	\leq \abs{ \{ A \subset G : d_1, d_2 \notin A - A\}  }.
\]
Putting everything together, we obtain the following lower bound
\begin{align}
\nonumber	\abs{\MSTD (G)} 
& \geq \abs{\{A \subset G : A - A \neq G, A + A = G \}}
\\ \nonumber&	= \abs{\bigcup_{d \in G'} D_d}
\\ \nonumber&	\geq \sum_{d \in G'} \abs{D_d} - \sum_{\substack{d_1, d_2 \in G' \\ d_1 \neq d_2}} \abs{D_{d_1} \cap D_{d_2}}
%\\ \nonumber&	\geq \sum_{d \in G'} \abs{ \{ A \subset G : d \notin A - A, A + A = G \}  } 
%\\ \nonumber &\qquad -	\sum_{\substack{d_1, d_2 \in G' \\ d_1 \neq d_2}} \abs{ \{ A \subset G : d_1, d_2 \notin A - A, A + A = G \}  } 
\\ \nonumber&	\geq \sum_{d \in G'} \(\abs{ \{ A \subset G : d \notin A - A \}  } 
					- \sum_{s \in G}  \abs{ \{ A \subset G : d \notin A - A, s \notin A + A \}  }  \)
\\ \label{eq:lower-union}&\qquad		- \sum_{\substack{d_1, d_2 \in G' \\ d_1 \neq d_2}} \abs{ \{ A \subset G : d_1, d_2 \notin A - A \}  }.
\end{align}
In the next section, we explain how to compute the individual terms on the RHS of \eqref{eq:upper-union} and \eqref{eq:lower-union}.

%%%%%%%%%%%%%%%%%%%%%%%%%%%%%%%%%%%%%%%%%%%%%%%%%%%%%%%%%%%%%%%%%%%%%%%%%%%%%%%%%%%%%%%%%%%%%%%

\section{Fibonacci indices of forbiddance graphs} \label{sec:forbiddance}

\subsection{Forbiddance graph} \label{sec:f-graph}

\begin{figure}[ht!] \centering
\begin{tikzpicture}[scale=1.5]
	\foreach \x in {0, 1, 2, ..., 7}
	{
		\node[P] (\x) at (\x *360/8 : 1) {};
		\node at (\x *360/8: 1.2) {$\x$};
	}
	\draw (0)--(1)--(2)--(3)--(4)--(5)--(6)--(7)--(0);
	\draw (1)--(3);	\draw (0)--(4);	\draw (5)--(7);
	\draw (2) edge[-,in = -45, out = -135,distance=10pt, loop] ();
	\draw (6) edge[-,in = 45, out = 135, distance=10pt, loop] ();
\end{tikzpicture}

\caption{The forbiddance graph of $\{A \subset \Zm 8 : 1 \notin A - A, 4 \notin A + A\}$. Note that the loops at $2$ and $6$ mean that neither vertex can be included in any independent set.\label{fig:Z8}}
\end{figure}
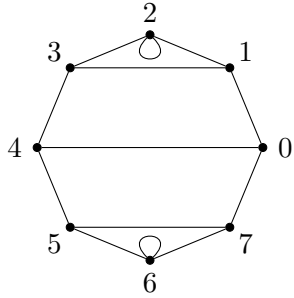

We would like to compute the cardinalities of collections of subsets of forms
\[
	\{ A \subset G : d \notin A - A \}, 
	\{ A \subset G : d \notin A - A, s \notin A + A \}, \text{ and } 
	\{ A \subset G : d_1, d_2 \notin A - A\}.
\]
More generally, for a collection of the form
\[
	\c A = \{ A \subset G : d_1, \dotsc, d_p \notin A - A; s_1, \dotsc, s_q \notin A + A \},
\]
we call $d_1, \dotsc, d_p$  the \emph{forbidden differences} and $s_1, \dotsc, s_q$ the \emph{forbidden sums}. We define the \emph{forbiddance graph} $\c G(\c A)$ to be the graph with $G$ as its sets of vertices, and an edge between two vertices (possibly coinciding) whenever their sum or difference is forbidden. In particular, if $2x$ is a forbidden sum, then there is a loop at $x$. See Figure \ref{fig:Z8} for an example of a forbiddance graph. If there are no forbidden sums, then the forbiddance graph is simply the undirected Cayley graph of the forbidden differences. The significance of forbiddance graphs is given in the following key observation. Recall that an independent set in a graph is a subset of the vertices with no two adjacent.

\begin{observation}
The elements of $\c A$ are in natural bijective correspondence with independent sets of vertices in the forbiddance graph $\c G(\c A)$.
\end{observation}

The number of independent sets in a graph $\sf G$ is known as the \emph{Fibonacci index} of $\sf G$, denoted $i({\sf G})$. The above observation implies that
\[
	\abs{\c A} = i(\c G (\c A)).
\]
Fibonacci indices were introduced by Prodinger and Tichy \cite{PT}. In chemistry, $i({\sf G})$ is known as the Merrifield-Simmons index \cite{MS}. Next we state some results about the Fibonacci index that we will use in the sequel.

\subsection{Fibonacci index} \label{sec:fib}

We use $F_n$ to denote the $n$-th Fibonacci number ($F_1=F_2 = 1$, $F_{n+2} = F_{n+1} + F_n$) and $L_n$ the $n$-th Lucas number ($L_1 = 1, L_2 = 3, L_{n+2} = L_{n+1} + L_n$). The path graph with $n$ vertices is denoted by $P_n$, the cycle graph with $n$ vertices is denoted by $C_n$, and the Cartesian graph product is denoted by $\square$.

The Fibonacci indices of many families of graphs are known (e.g., \cite{MW:indep-poly}). Here we state some results that are relevant to us. The derivation of the formulas can be found in the Appendix. (Assume $n \geq 1$; $C_1$ is the graph with one vertex and a loop.)
\begin{itemize}
	\item Path: $i(P_n) = F_{n+2}$.			
	\item Cycle: $i(C_n) = L_n = \vphi^n + (-\vphi)^{-n}$.			
	\item Ladder: $i(P_n \sq P_2) = \frac 1 2 \( (1 + \sqrt 2)^{n+1} + (1 - \sqrt 2)^{n+1} \)$.	 		
	\item Prism: $i(C_n \sq P_2) = (1 + \sqrt 2)^{n} + (1 - \sqrt 2)^{n} + (- 1)^n$.
	\item The Fibonacci index of a graph equals the product of the Fibonacci indices of its connected components.
\end{itemize}

We will also need an upper bound to the number of independent sets in a $d$-regular graph, which are graphs where vertex has degree $d$. This bound is provided by the following theorem which was conjectured implicitly by Alon \cite{Alon91}, explicitly by Kahn \cite{Kahn}, and recently proved by the author \cite{Zhao:indep}. See also \cite{GZ} for a different proof of Theorem \ref{thm:indep} when $d \leq 5$. We will only use the theorem for $d = 3$ and $4$. Note that we do not allow loops or multiple edges here.

\begin{theorem} \label{thm:indep}
For any $N$-vertex, $d$-regular simple graph $\sf G$
\[
	i({\sf G}) \leq (2^{d+1} - 1)^{N / (2d)}.
\]
\end{theorem}

\subsection{Forbidding one difference} \label{sec:1d}

We would like to determine the size of the collection
\[
	\c A = \{ A \subset G : d \notin A - A \}.
\]
The forbiddance graph $\c G(\c A)$ is easy to describe. Recall that $L_n$ is the $n$-th Lucas number.

\begin{lemma} \label{lem:cycles}
Let $d \in G$ be a nonzero element of order $\ell$. Let $\c A = \{ A \subset G : d \notin A - A \}$. Then the forbiddance graph $\c G(\c A)$ is a disjoint union of $\ell$-cycles $C_\ell$, and $\abs{\c A} = L_\ell^{\abs{G}/\ell}$.
\end{lemma}

\begin{proof}
The first statement is clear. The second statement follows from the facts listed in Section \ref{sec:fib}. 
\end{proof}

Let $\ord(d)$ denote the order of $d$ in $G$. The next result follows from the above lemma. Recall that $G'$ was defined in Section~\ref{sec:union-bound}.

\begin{lemma} \label{lem:cycles-sum}
We have
\[
	\sum_{d \in G'} \abs{  \{ A \subset G : d \notin A - A \} } 
					= \sum_{d \in G'} L_{\ord(d)}^{\abs{G} / \ord(d)}.
\]
\end{lemma}

The next lemma is an easy fact about Lucas numbers. 

\begin{lemma} \label{lem:L}
The sequence $\(L_{2n}^{1/2n}\)_{n\geq 1}$ is decreasing and the sequence $\(L_{2n-1}^{1/(2n-1)}\)_{n \geq 1}$ is increasing. Both sequences approach the limit $\vphi$. In particular, $L_2^{1/2} > L_4^{1/4} \geq L_n^{1/n}$ for all $n > 2$.
\end{lemma}

\begin{proof}
Recall that $L_n = \vphi^n + (-\vphi)^{-n}$, where $\vphi = \frac{1 + \sqrt 5}2$. So
\[
	L_{2n}^{1/2n} = (\vphi^{2n} + \vphi^{-2n})^{1/2n} = \vphi (1 + \vphi^{-4n})^{1/2n} \searrow \vphi,
\]
and
\[
	L_{2n-1}^{1/(2n-1)} = (\vphi^{2n-1} - \vphi^{-(2n-1)})^{1/(2n-1)} = \vphi (1 - \vphi^{-2(2n-1)})^{1/(2n-1)} \nearrow \vphi. \qedhere
\]
\end{proof}

The above fact is central to the dichotomy between the even case and the odd case in the theorem. In the sequel, we will show that, under the hypotheses of Theorem \ref{thm:main}, as $\abs{G} \to \infty$,
\begin{equation} \label{eq:union-asymp}
	\abs{\MSTD(G)} \sim \sum_{d \in G'} \abs{  \{ A \subset G : d \notin A - A \} } 
					= \sum_{d \in G'} L_{\ord(d)}^{\abs{G} / \ord(d)} \ ,
\end{equation}
In the case of even groups, the terms with $\ord (d) = 2$ dominate the RHS sum. For odd groups, every summand can be approximated from above by $\vphi^{\abs{G}}$, and the error is significant only when $\ord(d)$ is small. In the remainder of this section, we will analyze the other terms in \eqref{eq:lower-union} and show that they are insignificant compared to the RHS sum in \eqref{eq:union-asymp}.

\subsection{Other forbiddance graphs: even case} \label{sec:forb-even}

In this section we consider the case when $\abs{G}$ is even. We will show that the $L_2^{\abs{G}/2} = 3^{\abs{G}/2}$ term in \eqref{eq:union-asymp} asymptotically dominates the other terms in the sum in \eqref{eq:lower-union}.

\begin{lemma} \label{lem:even0}
If $d \in G$ has order greater than 2, then $\abs{\{ A \subset G : d \notin A - A \}} \leq 7^{\abs{G} / 4}$.
\end{lemma}

\begin{proof}
This follows from Lemmas \ref{lem:cycles} and \ref{lem:L}. Recall that $L_4 = 7$.
\end{proof}

Now let us consider the case when we forbid one sum and one difference. The structure of the forbiddance graph is given by the following lemma.

\begin{lemma} \label{lem:even1}
Let $d, s \in G$, where $d$ has order 2. Let $k$ denote the number of elements of $G$ of order 2.  Let $\c A = \{ A \subset G : d \notin A - A, s \notin A + A \}$. Then, after removing vertices with loops, the forbiddance graph $\c G (\c A)$ is a disjoint union of 4-cycles and at most $(k+1)/2$ copies of $P_2$.
\end{lemma}

\begin{proof}
Every connected component of $\c G (\c A)$ consists of the elements $\{x, x+d, s-x-d, s-x\}$, for $x \in G$, and possibly with some repeats. If all four elements are distinct, then this component is a 4-cycle. If $x = s-x$, then also $x+d = s-x-d$ (recall that $d$ has order 2), then this component consists of two connected vertices both with loops. If $x = s-x-d$, then also $x+d = s-x$, so the connected component is isomorphic to $P_2$. The number of $x \in G$ satisfying $2x = s-d$ is at most $k+1$, and hence at most $(k+1)/2$ components can be isomorphic to $P_2$.
\end{proof}

\begin{lemma} \label{lem:even2}
Let $d, s \in G$, where $d \neq 0$. Let $k$ denote the number of elements of $G$ of order 2. Then
\[ 
	\abs{\{ A \subset G : d \notin A - A, s \notin A + A \}} \leq
	\begin{cases}
		3^{(k+1)/2} \cdot 7^{\abs{G}/4}, & \text{if } d \text{ has order 2,} \\
		7^{\abs{G}/4}, & \text{if } d \text{ has order greater than 2.} \\
	\end{cases}
\]
\end{lemma}

\begin{proof}
In Lemma \ref{lem:even1}, note that the number of 4-cycles cannot exceed $\abs{G}/4$. Also recall that  $i(P_2) = 3$, $i(C_4) = 7$. The first case then follows immediately.

The second case follows from Lemma \ref{lem:even0} since
\[
	\abs{\{ A \subset G : d \notin A - A, s \notin A + A \}} 
	\leq \abs{\{ A \subset G : d \notin A - A\}}
	\leq 7^{\abs{G} / 4}. \qedhere
\]
\end{proof}

Next we move onto the case when we forbid two differences.

\begin{lemma} \label{lem:even3}
Let $d_1, d_2 \in G$ be distinct elements with order 2. Then $\{ A \subset G : d_1, d_2 \notin A - A\}$ has exactly $7^{\abs{G}/4}$ elements.
\end{lemma}

\begin{proof}
The connected components of the forbiddance graph are 4-cycles with vertices $x, x+d_1, x+d_1+d_2, x+d_2$ for $x \in G$. Note that all four elements are distinct. Since $i(C_4) = 7$, we see that the Fibonacci index of the forbiddance graph is $7^{\abs{G}/4}$.
\end{proof}

\begin{lemma} \label{lem:even4}
Assume that $\abs{G}$ is even. Let $d_1, d_2 \in G$ be distinct nonzero elements. Then 
\[
	\abs{ \{ A \subset G : d_1, d_2 \notin A - A\} } \leq 7^{\abs{G}/4}.
\]
\end{lemma}

\begin{proof}
If both $d_1$ and $d_2$ have order 2, then this follows from Lemma \ref{lem:even3}. Otherwise, suppose without loss of generality that $d_1$ has order greater than $2$, then from Lemma \ref{lem:even0} we get
\[
	\abs{ \{ A \subset G : d_1, d_2 \notin A - A\} }
	\leq \abs{ \{ A \subset G : d_1 \notin A - A\}} \leq 7^{\abs{G}/4}. \qedhere
\]
\end{proof}

\subsection{Other forbiddance graphs: odd case} \label{sec:forb-odd}

In this section we consider the case when $\abs{G}$ is odd. We will show that the RHS sum in \eqref{eq:union-asymp} dominates the RHS sum in \eqref{eq:lower-union}.

First we analyze the case when we forbid one sum and one difference.

\begin{lemma} \label{lem:prism-ladder}
Let $d \in G$ be of odd order $\ell > 1$, and $s \in G$. Let $\c A = \{A \subset G : d \notin A - A, s \notin A + A \}$. Then each connected component of the forbiddance graph $\c G(\c A)$, after removing vertices with loops, is either a  prism $C_\ell \sq P_2$ or a ladder $P_{(\ell-1)/2} \sq P_2$. Furthermore, if $n_P$ is the number of prism components and $n_L$ the number of ladder components, then $2n_P + n_L = \abs{G}/\ell$.
\end{lemma}

\begin{proof}
Let us start with the forbiddance graph of $\{ A \subset G : d \notin A - A\}$, which consists of $\abs{G}/\ell$ disjoint $\ell$-cycles, and then add edges of the form $(x, s-x)$, $x \in G$, to obtain $\c G(\c A)$. Let $x \in G$. There are two cases to consider.

\begin{description}
\item[Case 1.] $s - x$ does not belong to the same cycle as $x$. Then the connected component of $x$ in $\c G(\c A)$ is a prism $C_\ell \sq P_2$. The elements are connected as illustrated in Figure \ref{fig:prism}. %Otherwise, the connected component of $x$ is the 4-cycle $\{x, x+d, s-x-d, s - x\}$.
\item[Case 2.] $s -x $ belongs to the same cycle as $x$. Then the connected component of $x$ in $\c G(\c A)$ is the $\ell$-cycle with some ``parallel'' chords of the form $(x + jd, s - x - jd)$ added in. See Figure \ref{fig:chord}. Since we start with an odd cycle, there is exactly one vertex with a loop, and removing this vertex gives us the ladder $P_{(\ell-1)/2} \sq P_2$. 
\end{description}

\tikzstyle{P} = [draw, circle, black, fill, inner sep = 0pt, minimum width = 3pt]
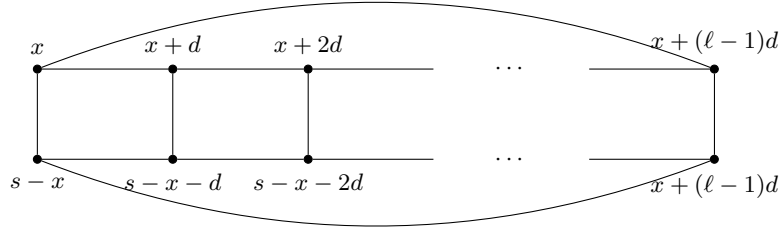
\begin{figure}[ht] \centering
\begin{tikzpicture}[xscale=1.8, yscale=1.2, font=\footnotesize]
	\node[P, label=above:$x$] (01) at (0,1) {};
	\node[P, label=above:$x+d$] (11) at (1,1) {};
	\node[P, label=above:$x+2d$] (21) at (2,1) {};
	\node (31) at (3,1){};
	\node at (3.5, 1) {$\cdots$};
	\node (41) at (4,1) {};
	\node[P, label=above:$x+(\ell-1)d$] (51) at (5,1) {};
	
	\node[P, label=below:$s-x$] (00) at (0,0) {};
	\node[P, label=below:$s-x-d$] (10) at (1,0) {};
	\node[P, label=below:$s-x-2d$] (20) at (2,0) {};
	\node (30) at (3,0){};
	\node at (3.5, 0) {$\cdots$};
	\node (40) at (4,0) {};
	\node[P, label=below:$x+(\ell-1)d$] (50) at (5,0) {};
	
	\draw (01)--(11)--(21)--(31);\draw (00)--(10)--(20)--(30);
	\draw (00)--(01);\draw (10)--(11);\draw (40)--(50);\draw (41)--(51);\draw (20)--(21);\draw (50)--(51);
	
	\draw (51) edge[bend right] (01);\draw (50) edge[bend left] (00);	
\end{tikzpicture}

\caption{The connected component of $x$ is the prism $C_\ell \sq P_2$ in Case 1 of the proof of Lemma \ref{lem:prism-ladder}. \label{fig:prism}}
\end{figure}

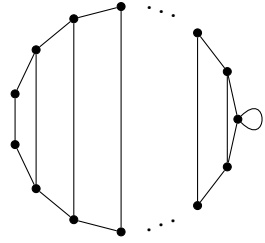
\begin{figure}[ht] \centering
\begin{tikzpicture} [scale=1.5]
	
	\node[P] (0) at (0:1) {};
	\node[P] (1) at (25:1) {};
	\node[P] (2) at (50:1) {};
	
	\node at (70:1) {\rotatebox{-20}{$\cdots$}};
	
%	\node[P] (3) at ( 92.5:1) {};
%	\node[P] (4) at (117.5:1) {};
%	\node[P] (5) at (142.5:1) {};
%	\node[P] (6) at (167.5:1) {};
%	\node[P] (7) at (-167.5:1) {};
%	\node[P] (8) at (-142.5:1) {};
%	\node[P] (9) at (-117.5:1) {};
%	\node[P] (10) at (-92.5:1) {};

	\node[P] (3) at ( 92:1) {};
	\node[P] (4) at (117:1) {};
	\node[P] (5) at (142:1) {};
	\node[P] (6) at (167:1) {};
	\node[P] (7) at (-167:1) {};
	\node[P] (8) at (-142:1) {};
	\node[P] (9) at (-117:1) {};
	\node[P] (10) at (-92:1) {};
	
	\node at (-70:1) {\rotatebox{20}{$\cdots$}};
	
	\node[P] (11) at (-50:1) {};
	\node[P] (12) at (-25:1) {};
	
	\draw (11)--(12)--(0)--(1)--(2);
	\draw (3)--(4)--(5)--(6)--(7)--(8)--(9)--(10);
	\draw (1)--(12);
	\draw (2)--(11);
	\draw (3)--(10);
	\draw (4)--(9);
	\draw (5)--(8);
	
	\draw (0) edge[-,in = 45, out = -45, distance=10pt, loop] ();
\end{tikzpicture}

\caption{The connected component of $x$ in Case 2 of the proof of the proof of Lemma \ref{lem:prism-ladder}.\label{fig:chord}}
\end{figure}

Therefore, after removing the vertices with loops, all connected components are prisms or ladders, thereby establishing the first claim. In the first case, two cycles combine to form a prism, and in the second case, a cycle transforms into a ladder. Since we start with $\abs{G}/\ell$ cycles, the claim $2n_P + n_L = \abs{G}/\ell$ follows.
\end{proof}

Using the notation from the lemma, we see that
\begin{align}
\nonumber	\abs{ \{A \subset G : d \notin A - A, s \notin A + A \}}
&	= i(P_{(\ell-1)/2} \sq P_2)^{n_L} \cdot i(C_\ell \sq P_2)^{n_P}
\\& 	\leq \max \left\{ i(P_{(\ell-1)/2} \sq P_2)^{\abs{G}/\ell},  i(C_\ell \sq P_2)^{\abs{G}/(2\ell)} \right\}. \label{eq:odd-ineq}
\end{align}
The following lemmas show that the above quantity is insignificant.

\begin{lemma} \label{lem:ladder-bound}
Let $\ell \geq 3$ be an odd integer. Then $i(P_{(\ell-1)/2} \sq P_2)^{1/\ell} < \sqrt{1 + \sqrt{2}}$.
\end{lemma}

\begin{proof}
Using results from Section~\ref{sec:fib} and routine algebraic manipulation, we obtain
\[
	i(P_{(\ell-1)/2} \sq P_2) = \frac 1 2 \( (1 + \sqrt 2)^{(\ell+1)/2} + (1 - \sqrt 2)^{(\ell+1)/2} \)
		< (1 + \sqrt 2)^{\ell / 2}. \qedhere
\]
%From Section \ref{sec:fib} we know that
%\[
%	i(P_n \sq P_2) = \frac 1 2 \( (1 + \sqrt 2)^{n+1} + (1 - \sqrt 2)^{n+1} \).
%\]
%So we need to prove that
%\[
%	\frac 1 2 \( (1 + \sqrt 2)^{n+1} + (1 - \sqrt 2)^{n+1} \) < (1 + \sqrt 2)^{n+\frac 1 2},
%\]
%which is equivalent to
%\[
%	(1 - \sqrt2)^{n+1} < (1 + \sqrt2)^{n} \( 2 \sqrt{1 + \sqrt{2}} - 1 - \sqrt 2 \).
%\]
%The last inequality is true because
%\[
%	(1 - \sqrt2)^{n+1} < \sqrt 2 - 1 < 2 \sqrt{1 + \sqrt 2} - 1 - \sqrt 2 < (1 + \sqrt2)^{n} \( 2 \sqrt{1 + \sqrt{2}} - 1 - \sqrt 2 \).
%	\qedhere
%\]
\end{proof}

\begin{lemma} \label{lem:prism-bound}
Let $\ell \geq 3$, then $i(C_\ell \sq P_2)^{1/(2\ell)} < 15^{1/6}$.
\end{lemma}

\begin{proof}
Since $C_\ell \sq P_2$ is a 3-regular simple graph, the lemma follows directly from Theorem~\ref{thm:indep}.
\end{proof}

We could have directly used the formula for $i(C_\ell \sq P_2)$ in Section \ref{sec:fib} to give a sharper bound, but the proof would be longer and the final conclusion would not be changed.

Using \eqref{eq:odd-ineq}, Lemmas \ref{lem:ladder-bound} and \ref{lem:prism-bound}, and noting that $\sqrt{1 + \sqrt 2} < 15^{1/6}$ we obtain the following upper bound for the case of forbidding one sum and one difference.

\begin{lemma} \label{lem:odd1d1s}
Suppose that $\abs{G}$ is odd, and $d, s \in G$ with $d \neq 0$, then
\[
	\abs{ \{A \subset G : d \notin A - A, s \notin A + A \}} < 15^{\abs{G}/6}.
\]
\end{lemma}

Next we consider the case of forbidding two differences.

\begin{lemma} \label{lem:odd2d}
Let $d_1, d_2 \in G$ be two nonzero elements, such that $2d_1 \neq 0$, $2d_2 \neq 0$, $d_1 \neq d_2$, $d_1 \neq - d_2$, then
\[
	\abs{ \{A \subset G : d_1, d_2 \notin A - A \} } \leq 31^{\abs{G}/8}.
\]
\end{lemma}

\begin{proof}
The forbiddance graph is simple and $4$-regular, since each vertex $x$ is adjacent to four distinct vertices: $x+d_1, x-d_1, x+d_2, x-d_2$. The result then follows from Theorem~\ref{thm:indep}.
\end{proof}

%%%%%%%%%%%%%%%%%%%%%%%%%%%%%%%%%%%%%%%%%%%%%%%%%%%%%%%%%%%%%%%%%%%

\section{Bounding the number of MSTD sets} \label{sec:bounding}

In this section we prove Theorem \ref{thm:main}. We use \eqref{eq:upper-union} and \eqref{eq:lower-union} to bound the number of MSTD sets in $G$ and use the results from the previous section to bound the individual terms. The even case and the odd case are analyzed separately.

\subsection{Even case}

Assume that $\abs{G}$ is even. Let $k$ denote the number of elements of order 2 in $G$. Then from \eqref{eq:upper-union} and Lemmas \ref{lem:cycles-sum} and \ref{lem:L}, we have
\begin{align}
	\nonumber \abs{\MSTD(G)} &\leq \sum_{d \in G'} \abs{ \{ A \subset G : d \notin A - A \} } 
	\\	\nonumber	&= \sum_{d \in G'} L_{\ord(d)}^{\abs{G}/\ord(d)}
	\\	\nonumber	&\leq k \cdot 3^{\abs{G}/2} + \abs{G} \cdot 7^{\abs{G}/4}
	\\	&= k \cdot 3^{\abs{G}/2} \(1 + \frac{\abs{G}}{k} \( \frac{7}{9} \)^{\abs{G}/4} \).	\label{eq:even-upper}
\end{align}

For the lower bound, we use some facts from Section \ref{sec:forb-even}. For any element $d \in G$ with order 2, Lemma \ref{lem:even2} gives
\begin{align*}
	&\abs{\{ A \subset G : d \notin A - A, A + A = G \} }
\\	& \quad\geq\abs{\{ A \subset G : d \notin A - A\}} - \sum_{s \in G} \abs{\{ A \subset G : d \notin A - A, s \notin A + A \}}
\\&\quad \geq 3^{\abs{G}/2} - \abs{G} 3^{(k+1)/2} \cdot 7^{\abs{G}/4}.
\end{align*}
From Lemma \ref{lem:even4} we know that for any distinct nonzero $d_1, d_2 \in G$,
\[
	\abs{\{ A \subset G : d_1, d_2 \notin A - A, A + A = G \}} \leq \abs{\{ A \subset G : d_1, d_2 \notin A - A\}} \leq 7^{\abs{G}/4}.
\]
Let $G_{(2)}$ be the subset of $G$ containing the elements of order 2, so that $\lvert G_{(2)} \rvert = k$. Then, from \eqref{eq:lower-union} we get
\begin{align}
\nonumber	\abs{\MSTD (G)} &\geq \sum_{d \in G'} \abs{\{ A \subset G : d \notin A - A, A + A = G \}} 
\\\nonumber	&\qquad  - \sum_{\substack{d_1, d_2 \in G' \\ d_1 \neq d_2}} \abs{\{ A \subset G : d_1, d_2 \notin A - A, A + A = G \}}
\\\nonumber	&\geq \sum_{d \in G_{(2)}} \abs{\{ A \subset G : d \notin A - A, A + A = G \} }
\\\nonumber	&\qquad		- \sum_{\substack{d_1, d_2 \in G' \\ d_1 \neq d_2}}  \abs{\{ A \subset G : d_1, d_2 \notin A - A, A + A = G\}}
\\\nonumber &\geq k\( 3^{\abs{G}/2} - \abs{G} 3^{(k+1)/2} \cdot 7^{\abs{G}/4} \) - \abs{G}^2 7^{\abs{G}/4}
\\\label{eq:even-lower} &= k \cdot 3^{\abs{G}/2} \(1 - \(\abs{G} \cdot 3^{(k+1)/2} + \frac{\abs{G}^2}{k}\) \( \frac{7}{9}\)^{\abs{G}/4} \).
\end{align}
By combining \eqref{eq:even-upper} and \eqref{eq:even-lower}, and using the notation of the even case of Theorem~\ref{thm:main}, we obtain that
\[
	 1 - \(\abs{G_n} \cdot 3^{(k_n+1)/2} + \frac{\abs{G_n}^2}{k_n}\) \( \frac{7}{9}\)^{\abs{G_n}/4} 
	\leq \frac{\abs{\MSTD (G_n)}}{k_n \cdot 3^{\abs{G_n}/2}}
	\leq 1 + \frac{\abs{G_n}}{k_n} \( \frac{7}{9} \)^{\abs{G_n}/4}.
\]
If $\abs{G_n} \to \infty$ and $\limsup_{n \to \infty} \frac{k_n}{\abs{G_n}} < 1 - \frac{1}{2}\log_3 7$, then letting $n \to \infty$ gives us 
\[
	\lim_{n\to \infty} \frac{\abs{\MSTD (G_n)}}{k_n \cdot 3^{\abs{G_n}/2}} = 1,
\]
thereby proving the even case of Theorem \ref{thm:main}.

\subsection{Odd case}

Now assume that $\abs{G}$ is odd. We use \eqref{eq:upper-union} and Lemma \ref{lem:cycles} to obtain an upper bound for $\abs{\MSTD(G)}$,
\begin{equation} \label{eq:odd-upper}
	\abs{\MSTD(G)} \leq \sum_{d \in G'} \abs{A \subset G : d \notin A - A} = \sum_{d \in G'} L_{\ord(d)}^{\abs{G}/\ord(d)}.
\end{equation}
For the lower bound, we use \eqref{eq:lower-union} and Lemmas \ref{lem:cycles}, \ref{lem:odd1d1s}, and \ref{lem:odd2d} to get
\begin{align}
\nonumber	\abs{\MSTD (G)} 
		&	\geq \sum_{d \in G'} \abs{ \{ A \subset G : d \notin A - A \}  } 
					- \sum_{d \in G'} \sum_{s \in G}  \abs{ \{ A \subset G : d \notin A - A, s \notin A + A \}  } 
\\ \nonumber &\qquad	- \sum_{\substack{d_1, d_2 \in G' \\ d_1 \neq d_2}} \abs{ \{ A \subset G : d_1, d_2 \notin A - A \}  }.
\\ \label{eq:odd-lower}		& \geq \(\sum_{d \in G'} L_{\ord(d)}^{\abs{G}/\ord(d)}\) - \abs{G}^2 \cdot 15^{\abs{G}/6} - \abs{G}^2 \cdot 31^{\abs{G}/8}.
\end{align}
By Lemma \ref{lem:L} and the fact that $\ord(d)$ is odd for every $d \in G$, we have $L_{\ord(d)}^{\abs{G}/\ord(d)} \geq L_3^{\abs{G}/3} = 4^{\abs{G}/3}$. Since $4^{1/3} > 15^{1/6} > 31^{1/8}$, \eqref{eq:odd-upper} and \eqref{eq:odd-lower} imply that
\begin{align} 
	\abs{\MSTD(G)} 
&	\sim \sum_{d \in G'} L_{\ord(d)}^{\abs{G}/\ord(d)}  \nonumber
\\&	\sim \frac{1}{2} \sum_{d \in G} L_{\ord(d)}^{\abs{G}/\ord(d)} \nonumber 
\\&	= \frac12 \sum_{d \in G} \(\vphi^{\ord(d)}  + (-\vphi)^{-\ord(d)} \)^{\abs{G} / \ord(d)} \nonumber 
\\&	= \frac12 \vphi^{\abs{G}} \sum_{d \in G} \(1 - \vphi^{-2\ord(d)} \)^{\abs{G} / \ord(d)}, \label{eq:odd-sum}
\end{align}
where the second line follows from $\ord (d) = \ord(-d)$. The second step is not an equality due to the negligible $d=0$ term. The last step uses the assumption that $\abs{G}$ is odd, so that $\ord(d)$ is odd. It remains to determine the asymptotics for the RHS of \eqref{eq:odd-sum}. We will need the following lemma.

\begin{lemma} \label{lem:bernoulli}
Let $n \geq k > 0$, then
\[
	1 - \(1 - \vphi^{-2k} \)^{n/k} \leq \min \left\{1, \frac{n}{k}\vphi^{-2k} \right\}.
\]
\end{lemma}

\begin{proof}
The inequality $1 - \(1 - \vphi^{-2k} \)^{n/k} \leq 1$ is obvious. The inequality
\[
	1 - \(1 - \vphi^{-2k} \)^{n/k} \leq  \frac{n}{k}\vphi^{-2k}
\]
follows directly from Bernoulli's inequality, which states that $(1+ x)^a \geq 1 + ax$ whenever $x > -1$ and $a \geq 1$ (Bernoulli's inequality can proved by checking the first derivative in $x$).
\end{proof}

Let $\c E(G)$ denote the subset of $G$ consisting of elements whose orders are less than $\log_{\vphi} \abs{G}$. So the hypotheses of Theorem \ref{thm:main} imply that $\abs{\c E(G)} = o(\abs{G})$.

\begin{lemma}
Let $G$ be an abelian group of odd order. Then
\[
	\abs{G} - \abs{\c E(G)} - 1 \leq \sum_{d \in G} \(1 - \vphi^{-2\ord(d)} \)^{\abs{G}/\ord(d)} \leq \abs{G}.
\]
\end{lemma}

\begin{proof}
The RHS inequality is clear since each term in the sum is at most 1. For the LHS, we need to prove that
\[
	\sum_{d \in G} \(1 - \(1 - \vphi^{-2\ord(d)} \)^{\abs{G}/\ord(d)}\) \leq \abs{\c E(G)} + 1.
\]
We separate the terms with small order from those with large order. We have
\[
	\sum_{d \in \c E(G)} \(1 - \(1 - \vphi^{-2\ord(d)} \)^{\abs{G}/\ord(d)}\) \leq \sum_{d \in \c E(G)} 1 = \abs{\c E(G)},
\]
and from Lemma \ref{lem:bernoulli} we get
\begin{align*}
	\sum_{d \notin \c E(G)} \(1 - \(1 - \vphi^{-2\ord(d)} \)^{\abs{G}/\ord(d)}\)
&	\leq \sum_{d \notin \c E(G)} \frac{\abs{G}}{\ord(d)}\vphi^{-2\ord(d)}
\\&	\leq \sum_{d \notin \c E(G)} \frac{\abs{G}}{\log_\vphi \abs{G}}\vphi^{-2\log_\vphi \abs{G}}
\\&	= \(\abs{G} - \abs{\c E(G)}\) \cdot \frac{\abs{G}}{\abs{G}^2\log_\vphi \abs{G}}
\\&	\leq 1
\end{align*}
This completes the proof of the lemma.
\end{proof}

Therefore, as $\abs{G} \to \infty$, if $\abs{\c E(G)} = o(\abs{G})$, then
\[
	\abs{\MSTD(G)}
	\sim \frac12 \vphi^{\abs{G}} \sum_{d \in G} \(1 - \vphi^{-2\ord(d)} \)^{\abs{G} / \ord(d)}
	\sim \frac12 \vphi^{\abs{G}} \abs{G}.
\]
This completes the odd case of Theorem \ref{thm:main}.

Recall the \emph{rank} of a finite abelian group is the smallest number of cyclic groups whose direct product gives the group. The following proposition shows that Theorem \ref{thm:main} can be applied when $\{G_n\}$ has uniformly bounded rank (this is already clear in the even case).

\begin{proposition} \label{prop:rank}
Let $r$ and $t$ be positive integers. Suppose that the finite abelian group $G$ has rank at most $r$, then the number of elements of $G$ with order less than $t$ is less than $t^{r+1}$.
\end{proposition}

\begin{proof}
In general, the number of elements of $\Zm a$ with order dividing $m$ is $\gcd(a, m)$. If an element of $\Zm {a_1} \x \cdots \x \Zm {a_r}$ has order dividing $m$, then each component has order dividing $m$, and the number of such elements is $\prod_i \gcd(a_i, m) \leq m^r$. Therefore, the number of elements of order dividing $m$ in $G$ is at most $m^r$. Summing over all $m < t$, we find that the number of elements with order less than $t$ is at most $1^r + 2^r + \cdots + (t-1)^r < t^{r+1}$.
\end{proof}

If the rank of $\{G_n\}$ is uniformly bounded by $r$, then $\abs{\c E(G_n)} < (\log_\vphi \abs{G_n} )^{r+1} = o(\abs{G_n})$, and thus the hypotheses of Theorem \ref{thm:main} in the odd case are satisfied.

%%%%%%%%%%%%%%%%%%%%%%%%%%%%%%%%%%%%%%%%%%%%%%%%%%%%%%%%%%%%%%%%%%%%%%%%%%%%%%%%%%%%%%%%%%%%%%%%%%

\section{Discussion and further directions} \label{sec:discussion}

We can ask for the rate of convergence in the asymptotics in Theorem \ref{thm:main}. A closer look at the proofs indicate that, in the even case, the ratio of the two quantities converges to 1 exponentially fast with respect to the size of the group. The odd case does not converge as quickly---our proofs show that the rate of convergence is at most linear in the proportion of elements with order less than $\log_\vphi \abs{G_n}$. We should be able to obtain a more precise asymptotic result for the odd case by including the number of elements with ``small'' order as a parameter. However, in the interest of keeping this paper short, we did not attempt such analysis.

We offer some ideas for generalizations and variations. The author \cite{Zhao:limit} recently studied the number of subsets of $\{1, \dots, n\}$ which miss a particular number of sums and a particular number of differences. It would be nice to investigate similar questions for the case of finite abelian groups: given nonnegative integers $s$ and $d$, how many subsets $A \subset G$ satisfies $\abs{A+A} = \abs{G} - s$ and $\abs{A-A} = \abs{G} - d$? We also ask for constructions of families of MSTD subsets of $\Zm n$, similar to the constructions of families of MSTD subsets of $\{1, \dots, n\}$ in the integers given recently by Miller, Orosz, and Scheinerman \cite{MOS} and the author \cite{Zhao:bidirectional}. There could be generalizations to other linear forms such as $A+A-A$; see \cite[Sec.~4]{MOS} for work in this direction in the integers case. Finally, we can consider variations where we choose a random subset of $G$ based on some other probability model, for instance, so that the number of expected elements in the chosen subset is not $\abs{G}/2$ but $O(\sqrt{\abs{G}})$; see Hegarty and Miller \cite{HM} for this analysis in the case of integers.

\appendix

\section{Computing the Fibonacci indices of certain graphs} \label{app:fib-compute}

In this appendix we prove the formula stated in Section \ref{sec:fib} for the number of independent sets in the path, cycle, ladder, and prism graphs.

\subsection{Path graph $P_n$}

Let $v$ be the first node in the path. If we do not include $v$ in the independent set, then rest of the independent set can be chosen as any independent set of $P_{n-1}$. On the other hand, if we include $v$ in the independent set, then the neighbor of $v$ cannot be included, and the rest of the independent set can be chosen as any independent set of $P_{n-2}$. Thus $i(P_{n}) = i(P_{n-1}) + i(P_{n-2})$. Checking the initial values, we find that $i(P_n) = F_{n+2}$, where $F_n$ is the $n$-th Fibonacci number ($F_1=F_2 = 1$, $F_{n+2} = F_{n+1} + F_n$). 

\subsection{Cycle graph $C_n$}

An argument analogous to the one above shows that
\[
	i(C_n) = i(P_{n-1}) + i(P_{n-3}) = F_{n+1} + F_{n-1} = L_n = \vphi^n + (-\vphi)^{-n}
\]
where $L_n$ is $n$-th Lucas number ($L_0 = 2, L_1 = 1, L_{n+2} = L_{n+1} + L_n$).

\subsection{Ladder graph $P_n \times P_2$}

We can compute $i(P_n \times P_2)$ using the transfer matrix method \cite[Sec.~4.7]{EC1}. The number of independent sets in $P_n \times P_2$ is equal to the number of walks of $n-1$ steps (starting at any vertex) in the following graph. Indeed, every such walk corresponds to a labeling of the vertices of $P_n \times P_2$ with $0$ and $1$ such that no two 1's are adjacent.

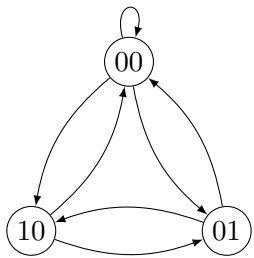
\begin{center}
\tikzstyle{V} = [draw, circle, inner sep = 2pt, minimum size = 15pt]
\begin{tikzpicture}[-latex, bend angle = 20, every loop/.style={-latex}]
	\node[V] (00) at (90 : 1.5) {$00$};
	\node[V] (01) at (-30 : 1.5) {$01$};
	\node[V] (10) at (210 : 1.5) {$10$};
	\draw (00) to[bend right] (10);
	\draw (00) to[bend right] (01);
	\draw (01) to[bend right] (10);
	\draw (01) to[bend right] (00);
	\draw (10) to[bend right] (00);
	\draw (10) to[bend right] (01);
	\draw (00) to[loop above] ();
\end{tikzpicture}
\end{center}

The transfer matrix (i.e., adjacency matrix) is
\[
	A = \begin{pmatrix} 1 & 1 & 1 \\ 1 & 0 & 1 \\ 1 & 1 & 0 \end{pmatrix}.
\]
The sum of all the entries of $A^n$ equals to the number of walks of $n$ steps in the above graph, which equals to $i(P_{n+1} \times P_2)$. We have
\[
	\sum_{n\geq 0} A^n x^n = (1 - Ax)^{-1} = \frac{1}{(1+x)(1 - 2x - x^2)} 
						\begin{pmatrix}
								1 - x^2 & x + x^2 & x + x^2 \\
								x + x^2 & 1 - x - x^2 & x \\
								x + x^2 & x & 1 - x - x^2  \end{pmatrix}.
\]
By taking the sum of all the entries of the matrix, we see that
\[
	\sum_{n \geq 0} i(P_{n+1} \times P_2) x^n 
	= \frac{3 + 4x + x^2}{(1+x)(1 - 2x - x^2)}  = \frac{3 + x}{1 - 2x - x^2}
	= \frac{1}{2x} \( \frac{1 + \sqrt 2}{1 - (1 + \sqrt 2) x} + \frac{1 - \sqrt 2}{1 - (1-\sqrt 2) x} - 2 \)
\]
By comparing the coefficient of $x^n$, we find that
\[
	i(P_n \x P_2) = \frac 1 2 \( (1 + \sqrt 2)^{n+1} + (1 - \sqrt 2)^{n+1} \).
\]

\subsection{Prism graph $C_n \times P_2$}

We can reuse the transfer matrix from the previous computation. Independent sets in $C_n \times P_2$ correspond to closed walks of length $n$, and so there are $\trace(A^n)$ of them. This can be computed as the number of sum of $n$-th powers of the eigenvalues of $A$. Hence,
\[
	i(C_n \times P_2) = (1 + \sqrt 2)^{n} + (1 - \sqrt 2)^{n} + (- 1)^n.
\]

%%%%%%%%%%%%%%%%%%%%%%%%%%%%%%%%%%%%%%%%%%%%%%%%%%%%%%%%%%%%%%%%%%%%%%%%%%%%%%%%%%%%%%%%%%%%%%%%%%

\section*{Acknowledgments}

This research was carried out at the University of Minnesota Duluth under the supervision of Joseph Gallian with the financial support of the National Science Foundation and the Department of Defense (grant number DMS 0754106), the National Security Agency (grant number H98230-06-1-0013), and the MIT Department of Mathematics. The author would like to thank Joseph Gallian for his encouragement and support. The author would also like to thank Nathan Kaplan for reading the paper and making valuable suggestions.
%The author would also like to thank Nathan Kaplan and Ricky Liu for reading the paper and making valuable suggestions.

%%%%%%%%%%%%%%%%%%%%%%%%%%%%%%%%%%%%%%%%%%%%%%%%%%%%%%%%%%%%%%%%%%%%%%%%%%%%%%%%%%%%%%%%%%%%%%%%

\bibliographystyle{amsplain}
\bibliography{../references}

\appendix 

\end{document}